\documentclass[a4paper,12pt,one side]{amsart}
\usepackage{amssymb,latexsym,xy,eucal,enumerate,tikz}
\usepackage{mathrsfs}
\usepackage{arydshln}
\usepackage{hyperref}
\usepackage{relsize}

\usepackage{pgfplots}
\usepackage{mathrsfs}
\usetikzlibrary{arrows.meta}
\usepackage{tikz}

\fontseries{m}
\usepackage{graphicx}
\graphicspath{ {./images/} }

\everymath{\displaystyle}

\theoremstyle{plain}
\newtheorem{lem}{Lemma}[section]
\newtheorem{lemma}[lem]{Lemma}
\newtheorem{theorem}[lem]{Theorem}

\newtheorem{prop}[lem]{Proposition}

\newtheorem{algorithm}[lem]{Algorithm}

\theoremstyle{definition}

\newtheorem{notation}[lem]{Notation}

\newtheorem{remark}[lem]{Remark}

\newtheorem{definition}[lem]{Definition}

\begin{document}
\title{ Enhanced power graph from the power graph of a group
}
\author{ Amar S. Pote}
\author{Ganesh S. Kadu}
\address{Amar S. Pote \newline \indent Sir Parashurambhau College, Savitribai Phule Pune University, Pune, India}
\email{amarpote0123@gmail.com}

\address{Ganesh S. Kadu  \newline \indent Department of Mathematics, Savitribai Phule Pune University, Pune, India}
\email{ganeshkadu@gmail.com}

\maketitle 
\begin{abstract} The power graph of a group $G$ is a graph with vertex set $G$, where two distinct vertices $a$ and $b$ are adjacent if one of $a$ and $b$ is a power of the other. Similarly, the enhanced power graph of $G$ is a graph with vertex set $G$, where two distinct vertices are adjacent if they belong to the same cyclic subgroup. In this paper we give a simple algorithm to construct the enhanced power graph from the power graph of a group  without the knowledge of the underlying group. This answers a question raised by Peter J. Cameron of constructing  enhanced power graph of group $G$ from its power graph. We do this by defining an arithmetical function on finite group $G$ that counts the number of closed twins of a given vertex in the power graph of a group. We compute this function and prove many of its properties.  One of the main ingredients of our proofs is the monotonicity of this arithmetical function on the poset of all cyclic subgroups of $G$.  
\end{abstract}

\noindent{ \small \textbf{Keywords}: {Group, Power graph, Enhanced power graph.}

\baselineskip 18truept
\section {Introduction}
The power graph $\mathscr{P}(G)$ of a group $G$ is a graph with vertex set  $G$, where two distinct vertices $a,b\in V(\mathscr{P}(G))$ are adjacent if and only if $\displaystyle \left < a\right > \subseteq \left <b \right >$ or  $\displaystyle \left < b \right > \subseteq \left <a\right>$. Power graphs were introduced by Kelarev and Quinn in \cite{Kelarev1, Kelarev2} as directed  graphs.  In \cite{Chakrabarty 1} Chakrabarty et al. introduced the undirected power graphs and proved that for a
finite group $ G$, $\mathscr{P}(G)$ is complete if and only if $G$ is  trivial or cyclic group of prime-power order. Cameron and Ghosh  in \cite{Cameron 3} showed that non-isomorphic finite groups may have isomorphic power graphs, but finite abelian groups with isomorphic power graphs must be isomorphic. \\ 
 \indent In \cite{Cameron 4} Aalipour et al. defined the enhanced power graph  of a group $G$. The enhanced power graph $\mathscr{P}_{e}(G)$ of a group $G$ is a graph having  vertex set as $G$, where two distinct vertices $a,b \in G$ are adjacent if there exists $c \in G $ such that $a,b \in  \left <c \right > $.  In \cite {Bera} Bera and Bhuniya proved that $\mathscr{P}_{e}(G)$ is complete if and only if $G$ is cyclic. 
 In \cite[Theorem 26]{Cameron 4}, it was  proved that if the power graphs of groups $G$ and $H$ are isomorphic, then their enhanced power graphs are also isomorphic, also see \cite{Zahi}.

 It is easy to see that, power graph of a group is contained in its enhanced power graph. Precursor to both the power graph and the enhanced power graph is the commuting graph $\mathscr C(G)$ of a group defined by Brauer and Fauler \cite{Brauer} in 1955 where two distinct vertices $a , b$ in $G$ are adjacent if they commute. Another related graph is the deep commuting graph  $\mathscr{D}_c(G)$ of a group $G$ where two  elements of $G$ 
 are joined  if and only if their inverse images in every central extension of $G$ commute.  One notes that the enhanced power graph of a group is a subgraph of its commuting graph. In fact, these graphs have the following graph hierarchy, $\mathscr{P}(G) \subseteq \mathscr{P}_e(G) \subseteq  \mathscr{D}_c(G) \subseteq \mathscr C(G).$ One of the central themes in this area has been to say when the two graphs in this graph  hierarchy are equal and, if they are not equal, then say something about the corresponding  difference graph, see \cite{Cameron 5}. For instance, in \cite[Theorem 28]{Cameron 4} it was proved that $\mathscr{P}(G)= \mathscr{P}_{e}(G)$ if and only if every cyclic subgroup of $G$ has prime power order. Such groups are known as EPPO groups and a complete list of EPPO groups is given in \cite[Theorem 1.7]{Cameron 6}. \\
  \indent In this paper, we consider the problem of constructing the enhanced power graph $\mathscr{P}_e(G)$ of a group $G$ from its power graph $\mathscr{P}(G)$  without the knowledge of the underlying group $G$. This question was raised by P. Cameron in \cite[Question 2]{Cameron 1}. We present a simple algorithm for constructing the enhanced power graph of $G$ from its power graph. The directed power graph can be constructed from the undirected power graph, see \cite{Bubboloni} and \cite{Das}. \\
   \indent To do this, we define  an arithmetical function on $G$ that counts the number of closed twins in $\mathscr{P}(G).$ We compute this function for finite groups and prove some of its properties. One of the main properties of this function is its monotonicity (Lemma \ref{monotonic}). We note that monotonicity of this arithmetical function  is crucial in proving the main theorems of this paper, namely,  Theorem \ref{main thm 1} and Theorem \ref{main thm 2}. Finally, using our main theorems we give a simple algorithm in Section $5$ for obtaining the power graph $\mathscr{P}(G)$ from the enhanced power graph $\mathscr{P}_e(G).$ As a result it is immediate that we have an algorithm for constructing difference graph directly from the power graph of group. In the last section we give some results on the difference graph.



\section{Preliminaries}
\noindent We begin by giving some notations which will be used throughout the paper.
\begin{notation} \label{notation} Let $G$ be a finite group. The following notations are used in this paper.
	\begin{enumerate}
		\item $U(\mathscr{P}(G))$ \hspace{0.37cm}: $\{ a \in \mathscr{P}(G) ~ \mid ~a \sim u$ for all $u \in \mathscr{P}(G)$\}. 
		\item $ N(a)$ \hspace{1.1cm}:  Neighbors of vertex $a$ in $\mathscr{P}(G)$  
		\item $\displaystyle \widetilde {N(a)}$\hspace{1.25cm}: $ \{b\in N(a)\mid N(a)\setminus\{b\}= N(b)\setminus\{a\}\}$ 
		\item $N_a$ \hspace{1.45cm}: $|\displaystyle \widetilde {N(a)}|+1$  
		\item $gen (\left <a \right >)$\hspace{0.63cm}: Set of generators of the cyclic subgroup $\left<a\right>$
	\end{enumerate}
Vertices in $\mathscr{P}(G)$ belonging to $\widetilde{N(a)}$ are called \textit{closed twins} of vertex $a$ in $\mathscr{P}(G)$.
\end{notation}

\begin{remark}
	For any group $G$, $V(\mathscr{P}(G))=V(\mathscr{P}_{e}(G))$ and $E(\mathscr{P}(G))\subseteq E(\mathscr{P}_{e}(G))$. 
\end{remark}

\begin{remark}
	Let $G$ be any finite group. Then,
	\begin{enumerate}
		\item Note that $ U(\mathscr{P}(G)) \neq \emptyset $,  as it contains identity element $\{e\}$, i.e., $\{e\}\subseteq  U(\mathscr{P}(G))$.
		\item If $G$ is cyclic, then all generators of $G$ are universal vertices,\\ i.e., $\{ gen (G)\} \subseteq U(\mathscr{P}(G))$.
	\end{enumerate}
\end{remark}
\noindent Following proposition by P. Cameron will be referred to frequently in this paper. 
\begin{prop} (\cite[Proposition 4]{Cameron 2})\label{Cam}
Let $G$ be a finite group with $|U(\mathscr{P}(G))|> 1.$ Then,
\end{prop}
\begin{enumerate}
	\item if $\mathscr P (G)$ is complete graph then, $G$ is cyclic of prime power order.
	\item if $\mathscr P (G)$ is incomplete then,  one of the following holds:
	\begin{enumerate}
	\item[\rm(i) ]	$G$ is cyclic of non-prime-power order $n\geqslant 6$, and $U(\mathscr{P}(G))$ contains the identity and generators of the $G$, so that $|U(\mathscr{P}(G))| = 1 + \phi(n)$.  
	\item[\rm(ii) ] $G$ is a generalized quaternion group of order $2^n$, and $U(\mathscr{P}(G))$ contains the identity and unique order $2$ element in $G$, so that	$|U(\mathscr{P}(G))|= 2.$
	 	\end{enumerate}
\end{enumerate} 

\begin{remark}
  In Proposition 2.1(2), G has to be the generalized quaternion group of order $2^n$. It cannot be any generalized quaternion group of order not equal to $2^n$ for any $n \geqslant 2$. This is because for generalized quaternion group of order not equal to $2^n$ one has, $|U(\mathscr{P}(G))|=1.$
\end{remark}
\noindent We now see that $G$ is non-cyclic  when $| U(\mathscr{P}(G))|= 1.$  
\begin{prop}
Let $G$ be a finite group. If $\mathscr{P}(G)$ is incomplete and $|U(\mathscr{P}(G))| = 1$, then $G$ is non-cyclic.
\end{prop}
\begin{proof}
On contrary suppose that $G$ is cyclic. By \cite[Theorem 2.12]{Chakrabarty 1}, , since $\mathscr{P}(G)$ is incomplete, we have $|G|\geqslant 6$. Since $G$ is cyclic, $U(\mathscr{P}(G))$ contains identity element $e$ and all generators of $G$. Since $|G|\geqslant 6$, $G$ has at least $2$ distinct generators. Therefore, $|U(\mathscr{P}(G))| \geqslant 3 $. This contradicts the assumption that $\mid U(\mathscr{P}(G)) \mid = 1$. Hence, $G$ is non-cyclic.
\end{proof}
\noindent The following is well known but we state it for the sake of completeness and for  referring.

\begin{prop}\label{cyclic}
Let $ H,K $ be cyclic subgroups of $G$ and suppose $H \subseteq K$. Then, 
\begin{enumerate}
	\item $|gen (H)|\leqslant |gen (K)|$.
	\item if $|gen (H)| = |gen (K)|$, then either $ |H| = |K|$ or $|K| = 2 |H|$. 
\end{enumerate}
\end{prop}
\begin{proof}
Suppose that $|H|=m $ and $|K|=n $, then $|gen (H)|= \phi (m)$ and $|gen (K)|=\phi (n)$.\newline Case (i): $m=n.$ Since $H \subseteq K$ we have, $H=K$. So, $|gen (H)|= |gen (K)|.$
  \newline Case (ii): $m\neq n$. Note $H \subsetneq K$ so, $m| n $, implying $n=mk$ for some $k \geqslant 2$. We claim that  $\phi (m)\leqslant \phi (n)$. Using Euler's totient function formula $\phi(n)=\phi(mk)= \phi(m) \phi(k)  \frac{d}{\phi(d)}$, where $d= \gcd (m,k)$.\newline Subcase (i): If $k=2$ then $\phi(k) =1 $ and $d= 1 $ or $2$, i.e., $\frac{d}{\phi(d)}= 1 $ or $2$. Thus either $\phi(n) = \phi (m)$ or $\phi(n) =  2 \phi (m)$. Therefore, $\phi (m)\leqslant \phi (n)$.
 \newline Subcase (ii): If $k > 2 $ then $\phi(k)\geqslant 2$. Since $ d \geq \phi (d)$, we have $\frac{d}{\phi(d)} \geq 1 $. Set $ r = \phi(k) \frac{d}{\phi(d)}.$ Note that $r \geqslant 2$. Thus, $\phi(n)= \phi(m)r$, with $r\geqslant 2$.
 So, $\phi (m) \leq \phi(n)$. This proves (1). \newline   To prove (2), we assume that $|gen (H)| = |gen (K)|$, i.e., $\phi(n)= \phi(m).$ Note that either $m=n$ or $ m \neq n.$  As in (1) above, we have two cases namely, $m=n$ or $m \neq n.$ If $m=n$ then we are done. If $m \neq n$ then by case (ii) above we have, $ m | n.$ Writing  $ n=mk$ as before we have $\phi(n)=\phi(mk)= \phi(m) \phi(k)  \frac{d}{\phi(d)}$, where $d= \gcd (m,k)$. As  $\phi(n)= \phi(m)$ we find that $\phi(k)  \frac{d}{\phi(d)}=1.$ Since $\phi(k) \geq 1$ and $\frac{d}{\phi(d)} \geq 1$ we get $\phi(k)=1.$ This gives $k=2,$ showing that $n=2m$ as required.
\end{proof}

\section{Function that counts closed twins and its properties}
\noindent We recall from the previous section that  $\displaystyle \widetilde {N(a)}= \big\{ b\in N(a)\mid N(a)\setminus\{b\}= N(b)\setminus\{a\}\big\} $ and $N_a =|\widetilde {N(a)}|+1. $ Note that $\displaystyle \widetilde {N(a)}$ is the set of closed twins of $a.$ This gives rise to an arithmetical function on $G$, namely,
\[ N: G \rightarrow \mathbb N_0 \]
\[   \quad \quad a \mapsto N_a\]
 We recall that the number $N_a$ counts the number of closed twins of vertex $a$ of the power graph. In this section, we compute this function for finite groups. We then prove some of the important properties of this function. In our main results, we utilize these properties to construct enhanced power graph from the power graph. We begin by showing that this counting function is lower bounded by Euler's totient function. In order to describe this we need to consider the Euler's totient function defined on any finite group $G$ as follows.
 \[\psi : G \rightarrow \mathbb N \]
 \[  \quad   \quad \quad \quad \quad \quad \quad \quad a \mapsto \psi(a):=\psi\big(|\left <a\right >|\big) \]
 Note that $\psi\big(|\left <a\right >|\big)= |gen (\left <a\right >)|.$ 
\begin{prop}\label{amar}
Let $G$ be a finite group and $a \in G.$ Then, 
 \[ gen (\left <a\right >)\setminus\{a\} \subseteq \displaystyle \widetilde {N(a)}.\]
 In particular, $  |gen (\left <a\right >)| \leq N_a $, i.e., $\psi(a) \leq N_a.$
\end{prop}
 \begin{proof}
 By definition, $\displaystyle \widetilde {N(a)}= \{b ~ | ~ ~ b \sim  a $ in $\mathscr{P}(G)$ and  $N(a)\setminus\{b\}=N(b)\setminus\{a\}\}$. Let $ b \in  gen (\left <a\right >)\setminus\{a\} $. Then, $ \left<a\right> = \left <b\right>$ and $b \sim a $ in $\mathscr{P}(G)$. We show that $ b \in  \widetilde{N(a)}$. To do this we need to show that   $N(a)\setminus\{b\}=N(b)\setminus\{a\}.$ So, let $t \in N(a)\setminus\{b\}.$  By definition of the power graph, $ \left<t\right> \subseteq \left <a\right>$ or $ \left<a\right> \subseteq \left <t\right>$. Since, $ \left<a\right> = \left <b\right>$, we have $\left<t\right> \subseteq \left <b\right>$ or $ \left<b\right> \subseteq \left <t\right>$. Thus, $t\sim b $ in $\mathscr{P}(G)$. This shows that $N(a)\setminus\{b\} \subseteq N(b)\setminus\{a\}.$ Similarly, $ N(b)\setminus\{a\} \subseteq  N(a)\setminus\{b\}.$ Hence, $N(a)\setminus\{b\}=N(b)\setminus\{a\}$ showing that $ b \in \displaystyle \widetilde {N(a)}$.  This completes the proof.          
 \end{proof}

\begin{prop}\label{non-prime-power}
Let $G$ be a finite non-cyclic group and $\left<a\right> \leqslant G$, where $a \neq e$ is an element of non-prime-power order. Then, $\displaystyle \widetilde {N(a)}= gen (\left <a\right >)\setminus\{a\}$ and hence, $ N_a= | gen (\left <a \right >)|=\psi(a)$.
\end{prop}
\begin{proof}
Let $a \neq e$ in $G$ be an element of non-prime-power order.  Since $a$ is an element of non-prime power order we have, $| \langle a \rangle|= p_1^{\alpha_1} p_2^{\alpha_2}... p_n^{\alpha_n}$, where $p_i$'s are all distinct primes and at least two $\alpha_i$'s are non-zero.  By definition, $N_a= |\widetilde {N(a)}| + 1$.  To prove $ N_a= | gen (\left <a \right >)|$, it is enough to show that  $\displaystyle \widetilde {N(a)}= gen (\left <a\right >)\setminus\{a\}.$ Already, by Proposition \ref{amar} we have, $ gen (\left <a\right >)\setminus\{a\} \subseteq \displaystyle \widetilde {N(a)}.$  We now show that $ \displaystyle \widetilde {N(a)} \subseteq  gen (\left <a\right >)\setminus\{a\}.$ So, let $b  \in \widetilde {N(a)}$ and suppose on the contrary that $ b \notin  gen (\left <a\right >)\setminus\{a\}.$    We first  show that element $b$ cannot be the identity $G.$ Then in this we have $e \in \widetilde{ N(a)}.$ This gives $a \in U(\mathscr{P}(G)).$ As $ a \neq e$ we have  $|U(\mathscr{P}(G))| > 1.$ Since $G$ is non-cyclic we have by \cite[Theorem 2.12]{Chakrabarty 1} $\mathscr P(G)$ is incomplete. Hence, by Proposition \ref{Cam}, we get that $G$ is a generalized quaternion group of order $2^n$ and that $a$ is unique element of order $2$ in $G.$ This contradicts the fact that $a$ is an element of non-prime-power order. Hence, $b \neq e.$ Since $ b\in \widetilde {N(a)}$, $b \sim a $ in $\mathscr{P}(G)$ by definition of $\widetilde {N(a)}$. Therefore, either $\left<b\right> \subseteq \left <a\right>$ or $ \left<a\right> \subseteq \left <b\right>$. Suppose first that $\left<b\right> \subseteq \left <a\right>.$  We now make two cases on $b$ according as order of $b$ is a prime power or a non-prime power.\\
\underline{Case (i)}: Suppose $|\langle b \rangle|= p_{i}^{\beta_{i}}$ for some $i$ and $\beta_{i}\leqslant {\alpha_i}$. (i.e., $b$ is an element of prime power order).
Since $\langle a \rangle $ is cyclic of order $p_1^{\alpha_1} p_2^{\alpha_2}... p_n^{\alpha_n}$, there exists $s \in \langle a \rangle $ such that  $|\langle s \rangle|= p_j^{\alpha_j}$ with $i \neq j$. Note that $s \neq b.$ Also, since $s \in \langle a \rangle $  we have  $s \sim a$ in $\mathscr{P}(G).$  Since orders of $b$ and $s$ do not divide each other it follows that $ b \nsim s $, i.e., $b$ and $s$ are non-adjacent in $\mathscr{P}(G).$ Thus, $s \in N(a) \setminus \{b\}$ but $ s \notin N(b) \setminus \{a\}$ showing that $ N(a) \setminus \{b\} \neq N(b) \setminus \{a\} $. This contradicts $b \in \widetilde {N(a)}.$\\
\noindent\underline{Case (ii)}: Suppose $|\langle b \rangle|\neq  p_{i}^{\gamma_{i}}$ for any prime $p_{i}$ and $\gamma_i \geq 1$.  
Since $\left<b\right> \subseteq \left <a\right>$ we have, $|\langle b \rangle|= p_1^{\beta_1} p_2^{\beta_2}... p_r^{\beta_r}$ for some $\beta_i \leq \alpha_i.$ Further, as $|\langle b \rangle|$ is an element of non-prime-power order, at least two $\beta_{i}$'s  are non-zero. Since $b \in \langle a \rangle $ and $ b \notin  gen (\left <a\right >)\setminus\{a\}$ there exists $\beta_j$ for some $j$ such that  $\beta_j < \alpha_j.$ We make two cases: either $r= n$ or $r \neq n$.\\    
\underline{Subcase (i)}: Suppose $r\neq n $. (i.e., $r < n $). Since $\langle a \rangle $ is cyclic there exists $s \in \langle a \rangle $  such that $|\langle s \rangle|= p_{r+1}^{\beta_{r+1}}$. Therefore, $s \sim a $ but $s \nsim b$ in $\mathscr{P}(G)$, contradicting $b \in \widetilde {N(a)}$.\\
\underline{Subcase (ii)}: Suppose $r = n $. As noted above as $ b \notin  gen (\left <a\right >)\setminus\{a\},$ there exists $\beta_j$ for some $j$ such that  $\beta_j < \alpha_j.$ So, there exists  $s \in \langle a \rangle $  such that $|\langle s \rangle|= p_j^{\alpha_j}$ with $\beta_j < \alpha_j$. Clearly  $ s\sim a.$ Further, since the orders of $s$ and $b$ do not divide each other we have, $s \nsim b$ in $\mathscr{P}(G)$.  This contradicts $b \in \widetilde {N(a)}$. \\
\indent Contradictions in both the cases above shows that $\left<b\right> \nsubseteq \left <a\right>.$  Therefore we are in the other case where $\left<a\right> \subseteq \left <b\right>.$ Note that since $a$ has non-prime power order implies that $b$ has non-prime power order and   $| \langle a \rangle|= p_1^{\alpha_1} p_2^{\alpha_2}... p_n^{\alpha_n}$ divides order of $| \langle b \rangle|.$ Thus, $| \langle b \rangle|=p_1^{\gamma_1} \ldots p_n^{\gamma_n}\ldots  p_m^{\gamma_m} $ where  $ \alpha_i \leq \gamma_i$ for all $1 \leq i \leq  n$ and $ m \geq n.$ We now make two cases, namely $m=n$ and $m > n.$\\
\underline{Case (i)} $m=n$:  Note that since $b \notin  gen (\left <a\right >)\setminus\{a\} $  we have that  $\langle a \rangle \subsetneq \langle b \rangle. $ So,  $|\langle a \rangle| < |\langle b \rangle|.$     Thus, there exists $j$ such that $  \alpha_j<\gamma_j.$ Note that there exists $s \in \langle b \rangle $ such that $|\langle s \rangle|= p_j^{\gamma_j}.$ Clearly $ s \sim b$.  Since the orders of $s$ and $a$ do not divide each other it follows that $s \nsim a.$  This contradicts $b \in \widetilde {N(a)}.$ \\
\underline{Case (ii)} $m > n$: In this case there exists $s \in \langle b \rangle $ such that $|\langle s \rangle|= p_m^{\gamma_m}.$ Note that $ s\sim b$ but $ s \nsim a$ as orders of $s$ and $a$ do not divide each other. This contradicts $b \in \widetilde {N(a)}.$ \\
\indent Contradictions in the both the cases above therefore show that $\left<a\right> \nsubseteq \left <b\right>.$ Thus, we have both  $\left<a\right> \nsubseteq \left <b\right>$ and  $\left<b \right> \nsubseteq \left <a \right>.$ This violates  $ b \notin \widetilde {N(a)}.$ Hence, there is a overall contradiction to the fact that 
$b \notin  gen (\left <a\right >)\setminus\{a\} $. This proves $b \in  gen (\left <a\right >)\setminus\{a\} $ showing that  $ \displaystyle \widetilde {N(a)} \subseteq  gen (\left <a\right >)\setminus\{a\}.$  Hence, $\displaystyle \widetilde {N(a)}= gen (\left <a\right >)\setminus\{a\}.$ Thus, $ N_a= | gen (\left <a \right >)|$ as desired.          
\end{proof}
 Next, we compute the twin counting function at the identity element of the group.  
 \begin{prop}\label{Ne}
 Let $G$ be a finite non-cyclic group. For identity $e$ of $G$ we have,
\[N_e = \left\{ \begin{array}{rcl}
	2 & \mbox{if G is a generalized quaternion of order $2^n$}\\
	 1 & \mbox{else.} 
\end{array}\right.\]
 \end{prop}

\begin{proof}
	It is easy to see that $\widetilde {N(e)} = U(\mathscr P(G)) \setminus \{e\}.$ Hence, $$N_e= |\widetilde {N(e)}| +1 = | U(\mathscr P(G)) \setminus \{e\}| +1 = | U(\mathscr P(G))|.$$  
	If $G$ is generalized quaternion group of order $2^n$ then by \cite[Theorem 2.12]{Chakrabarty 1} we have $N_e=| U(\mathscr P(G))| =2.$ Suppose now that $G$ is not a generalized quaternion group of order $2^n.$ We claim that $| U(\mathscr P(G))| =1.$  Suppose on the contrary that  $| U(\mathscr P(G))| > 1.$ Since $G$ is non-cyclic, it follows from  \cite[Theorem 2.12]{Chakrabarty 1} that $\mathscr P(G)$ is incomplete. By Proposition \ref{Cam} it follows that $G$ is a generalized quaternion group of order $2^n$ contradicting the assumption that $G$ is not a generalized quaternion group of order $2^n.$ Hence, $| U(\mathscr P(G))| =1$ in this case. This completes the proof.
	\end{proof}
\begin{lemma}
Let $G$ be a finite group and $\left<h\right>\leqslant G$, where $|\left<h\right>|=p^\alpha$, $p$ prime and $\alpha > 0$. Then, there does not exist $k\in \widetilde {N(h)}$ such that $|\left<k\right>|= p^i q^j t$, where $q \neq p$ is a prime and $i, j,t \geqslant 1.$  In other words, if $k\in \widetilde {N(h)}$ then $|\left<k\right>|= p^s$ for some $ s \geq 0.$  
\end{lemma}
\begin{proof}
Suppose on the contrary that there exists $k\in \widetilde {N(h)} $ and  primes $q \neq p$ such that $|\left<k\right>|= p^i q^j t$, where  $i, j,t \geqslant 1.$  Since $\left<k\right>$ is cyclic, there exists element $ l \in \left<k\right>$ of order $q $. Clearly, $l \sim k$ in $\mathscr P(G).$ Since the orders of $l$ and $h$ do not divide each other it follows that   $l \nsim h .$ Thus, $l$ is a neighbor of $k$ but $l$ is not a neighbor of $h.$ However, $k\in \widetilde {N(h)}$ implies $k$ and $h$ have same neighbors. This is a contradiction. Therefore, no such  $k\in \widetilde {N(h)}$ exists.
\end{proof}

\begin{prop}\label{Euler}
Let $G$ be a finite non-cyclic group and $h \in G$ be an element of prime-power order. Suppose there exists $ k \in G$ such that $\left<h\right>\subsetneq \left<k\right>$ and $k$ is  an element of non-prime-power order. Then, $\displaystyle \widetilde {N(h)}= gen (\left <h\right >)\setminus\{h\}  $ and hence, $N_h= | gen (\left <h \right >)| =\psi(h).$ 
\end{prop}
\begin{proof}
 Let $|\left<h\right>|= p^{\alpha}$	for some $\alpha \geq 1.$ Since  $\left<h\right>\subsetneq \left<k\right>$  and $k$ is an element of prime-power order we have,  $|\left<k\right>|= p^i q^{j} t$  for some $p \neq q$  distinct primes and  $ i, j,t\geqslant 1$.
To prove the proposition it is enough to prove $ gen (\left <h\right >)\setminus\{h\} = \displaystyle \widetilde {N(h)}.$ Already, from Proposition 2.6 we have, $ gen (\left <h\right >)\setminus\{h\} \subseteq \displaystyle \widetilde {N(h)}$. Let $l \in \widetilde {N(h)}.$ We show that $ l \in gen (\left <h\right >)\setminus\{h\}.$ Suppose on the contrary that $l \notin gen (\left <h\right >)\setminus\{h\}.$  Since  $l \in \widetilde {N(h)}$ by Lemma 2.9 we have, $|\left<l\right>|= p^s$ for some $ s \geq 0.$  Note that $l \sim h$ hence, either $ \left <h\right > \subseteq \left <l\right >$ or $\left <l\right > \subseteq \left <h\right >.$ Since $l \notin gen (\left <h\right >)\setminus\{h\}$ we have that $ \left <h\right > \subsetneq \left <l\right >$ or $\left <l\right > \subsetneq \left <h\right >.$ Thus,  $|\left<l\right>| \neq  p^\alpha =|\left<h\right>|,$ i.e., $s \neq \alpha.$ We make two cases, namely  $ \left <h\right > \subsetneq \left <l\right >$ or $\left <l\right > \subsetneq \left <h\right >.$ \\
 Case (i): Suppose $\langle l \rangle \subsetneq \langle h \rangle$.  In this case we have, $s < \alpha.$  If $s=0$, then $l =e$, where $e$ is the identity element of $G$. Since $e$ is adjacent to all elements in $V(\mathscr{P}(G))$ and $l=e \in \widetilde {N(h)}$ we have,  $h\in U(\mathscr{P}(G))$. So, $|U(\mathscr{P}(G))| >1.$ Note that since $G$ is non-cyclic, by  \cite[Theorem 2.12]{Chakrabarty 1}, $\mathscr{P}(G)$ is incomplete.  Since $G$ is non-cyclic and $\mathscr{P}(G)$ is incomplete,  by Proposition \ref{Cam}(2)(ii) $G$ is a generalized quaternion group of order $2^n.$ But $|\left<k\right>|= p^i q^{j}t$  where  $ i, j,t\geqslant 1$ with $p\neq q$ distinct primes which is clearly not possible in a group of order $2^n.$  Hence, $s\neq 0$. Since $\langle h \rangle \subsetneq \langle k \rangle$, we have $\langle l \rangle \subsetneq  \langle h \rangle \subsetneq \langle k \rangle$. Note that $\langle h \rangle \subsetneq \langle k \rangle$ we have, $\alpha \leq i.$ Thus, $s < \alpha \leq i.$ As $\langle k \rangle $ is cyclic, there exists $k' \in \langle k \rangle$  such that $|\langle k' \rangle| = p^s q$. Since $|\langle l \rangle|$ divides $|\langle k' \rangle|$ and $\langle l \rangle$,  $\langle k' \rangle$ are subgroups of the cyclic group $\langle k \rangle$ we have, $\langle l \rangle \subseteq \langle k' \rangle$ and  so $l \sim k'$. Further, since the orders of $h$ and $k'$ do not divide each other we obtain, $ k'\nsim h.$ Thus, $k'$ is a neighbor of $l$ but not a neighbor of $h,$  contradicting $l \in \widetilde {N(h)}$. \\ 
Case (ii): Suppose $\langle h \rangle \subsetneq \langle l \rangle$. In this case we have, $\alpha  < s$. As $\langle k \rangle $ is cyclic, there exists $k' \in \langle k \rangle$  such that $|\langle k' \rangle| = p^\alpha q$. Since $|\langle h \rangle|$ divides $|\langle k' \rangle|$ and $\langle h \rangle$,  $\langle k' \rangle$ are subgroups of the cyclic group $\langle k \rangle$ we have, $\langle h \rangle \subsetneq \langle k' \rangle$ and  so $h \sim k'$. Further, since the orders of $l$ and $k'$ do not divide each other we obtain, $ k'\nsim l.$ Thus, $k'$ is a neighbor of $h$ but not a neighbor of $l.$  This contradicts $l \in \widetilde {N(h)}$.

Contradictions in both the cases (i) and (ii) above shows that, $ l \in gen (\left <h\right >)\setminus\{h\}.$ This shows that   $ \displaystyle \widetilde {N(h)} \subseteq gen (\left <h\right >)\setminus\{h\} $ proving the desired equality. Hence, $N_h= | gen (\left <h \right >)|$.             
\end{proof}

We now summarize the compuation of the arithmetical function $ N : G \rightarrow \mathbb N$ which counts the number of closed twins of a vertex in $G.$
\begin{theorem}
	Let $G$ be a finite group of order $|G|=n$ and $ a\in G.$ Then the numbers $N_a$ are given as follows,
	\begin{enumerate}
		
		\item If $G$ is a cyclic group of prime-power order, then for each $a\in G$,  
		$N_a= n.$
		\item  If $G$ is cyclic group of non-prime-power order, then 
		\[N_a = \left\{ \begin{array}{rcl}
			\phi(n)+1, & \mbox{if $a$ is the identity or generator of $G$,}\\
			\phi(d), & \mbox{otherwise, where $d$ is order of $a.$} 
		\end{array}\right.\] 
		\item If $G$ is a non-cyclic group then,
		\begin{enumerate}
			
			\item For the identity $e\in G$,
			\[N_e = \left\{ \begin{array}{rcl}
				2 & \mbox{if G is  generalized quaternion of order $2^n$}\\
				1 & \mbox{otherwise.} 
			\end{array}\right.\]
			\item   $a \neq e$ is an element of non-prime-power order $d$, then $ N_a= \phi(d).$\\
			\item  $a \neq e$ has prime-power order $d$, and there exists an element of non-prime-power order  $ b \in G$ such that $\left<a\right>\subsetneq \left<b\right>$, then $ N_a= \phi(d).$ 
			\end{enumerate}
			\end{enumerate}
\end{theorem}

\begin{proof}
	(1) If $G$ is a cyclic group of prime-power order then it follows from \cite[Theorem 2.12]{Chakrabarty 1} see that $U(\mathscr{P}(G))=G.$ Thus, for any $a \in G$ we have, $\widetilde {N(a)}= G \setminus \{a\}.$ Hence, $N_a = |\widetilde {N(a)}|+1= |G| -1 +1= |G|.$\\
	(2) Assume $G$ is a cyclic group of non-prime-power order. Clearly $n=|G| \geq 6.$ By 	Proposition \ref{Cam}(2)(i), $U(\mathscr{P}(G)) = gen(G) \cup \{e\}$.  If $ a \in gen(G) \cup \{e\}$ then we have, $\widetilde {N(a)}= U(\mathscr{P}(G)) \setminus \{a\}.$ So, $N_a =|\widetilde {N(a)}|+1 =|U(\mathscr{P}(G)) \setminus \{a\}|+1=|U(\mathscr{P}(G))|-1+1=|U(\mathscr{P}(G))| = \phi(n)+1,$ where the last equality follows from Proposition \ref{Cam}(2)(i).\\ Next, suppose  that $a \notin gen(G) \cup \{e\}$ and order of $a$ is $d.$ There are two cases to consider.\\
	Case (i): Suppose $a$ has prime-power order. 
	Since $G$ is cyclic, $G=\left<b\right>$ with $b$ of non-prime-power order and $\left<a\right> \subsetneq \left<b\right>$.  Using arguments similar to those in the proof of Proposition \ref{Euler} we get, $N_a = |gen (\left <a \right >)|$, i.e., $N_a =\phi(d).$ \\ 
		Case(ii): Suppose $a$ has non-prime-power order. Using arguments similar to those in the proof of proposition \ref{non-prime-power} we get, $N_a = |gen (\left <a \right >)|$. i.e., $N_a =\phi(d)$.\\
		(3) Suppose $G$ is a non-cyclic group. In this case (a) follows from Proposition \ref{Ne}, (b) follows from Proposition \ref{non-prime-power} and (c) follows from Proposition \ref{Euler}.


\end{proof}

\begin{definition}\label{def max chain}
	Let $G$ be group with identity $e$. A  chain $\mathcal C: \{e\}=C_0 < C_1 <...< C_k$ of subgroups of $G$ is called a \textit{cyclic chain} if each $C_i$ is a cyclic subgroup of $G$  such that each inclusion in $\mathcal C$ is strict.  Such a chain is called \textit{maximal cyclic chain} of  subgroups if it is not properly contained in any longer chain of cyclic subgroups of $G$. Equivalently, for each $i$, there is no cyclic subgroup $H$ of $G$ such that $C_{i-1} < H < C_{i} $ for $i=1, 2, \ldots, k$ and $C_k$ is not properly contained in any cyclic subgroup of $G$.
\end{definition}

\begin{definition}
	A collection of subgroups of a group $G$ is called a \emph{chain} (or \emph{totally ordered by inclusion}) 
	if every two subgroups in the collection are comparable; that is, for any $H_1$ and $H_2$, either $H_1 \subseteq H_2$ or $H_2 \subseteq H_1$.
\end{definition}
\begin{lemma}\label{lab1}
	Let $G$ be a finite group and let $h \neq k \in  G$ such that $\left <h \right >= \left < k\right >$. Then, $\widetilde {N(h)}\setminus \{k\} = \widetilde {N(k)}\setminus \{h\}$ and hence, $N_h=N_k$. 
\end{lemma}
	\begin{proof}
		 Note first that $h \sim k$. It is enough to show that $\widetilde {N(h)}\setminus \{k\} = \widetilde {N(k)}\setminus \{h\}.$ To do this, we let $x \in \widetilde {N(h)}\setminus \{k\}$ and show that $x \in \widetilde {N(k)}\setminus \{h\}$, i.e., $x$ and $k$ have same neighbors.  Since $\left <h \right >= \left < k\right >$, we have $h \sim k $ in $\mathscr{P}(G)$, which implies $x \sim k.$ Suppose $t \sim x $, we need to show $ t\sim k$. Since $ t \sim x $ and $ x \in \widetilde {N(h)}\setminus \{k\}$, we have $ t \sim h$, i.e, $\displaystyle \left < t\right > \subseteq \left <h \right >$ or  $\displaystyle \left < h \right > \subseteq \left <t \right>$. But $\left <h \right >= \left < k\right >$, so  $\displaystyle \left < t\right > \subseteq \left <k \right >$ or  $\displaystyle \left < k \right > \subseteq \left <t \right>$. By definition, this means $ t\sim k.$ Thus, every neighbor of $x$ is also a neighbor of $k$. Conversely, if $t \sim k $, we need to show $ t\sim x $. Since $ t \sim k $, we have $\displaystyle \left < t\right > \subseteq \left <k \right >$ or  $\displaystyle \left < k \right > \subseteq \left <t \right>$. Using $\left <h \right >= \left < k\right >$, it follows that $\displaystyle \left < t\right > \subseteq \left <h \right >$ or  $\displaystyle \left < h \right > \subseteq \left <t \right>$, which implies $ t\sim h $. Since $ x \in \widetilde {N(h)}\setminus \{k\}$ we have,  $ t \sim x $. Hence, we conclude that $\widetilde {N(h)}\setminus \{k\} \subset \widetilde {N(k)}\setminus \{h\}.$  By symmetry, $\widetilde {N(k)}\setminus \{h\} \subset \widetilde {N(h)}\setminus \{k\}$. So, $\widetilde {N(h)}\setminus \{k\} = \widetilde {N(k)}\setminus \{h\}.$ This completes the proof.
	\end{proof}

\begin{lemma}\label{lab2}
	Let $G$ be any finite group. Let $\left<h\right>, \left<k\right>\leqslant G$ with $\left<h\right>\subsetneq  \left<k\right>$ and let $ x \in \widetilde {N(h)}.$ Then, there exists a cyclic chain $\mathcal C$ containing $x, h$ and $k.$
	\begin{proof}
		 Suppose  $\left<h\right>\subsetneq \left<k\right>$ and $x \in \widetilde {N(h)}.$  It is clear that $h \sim k $ in $\mathscr{P}(G)$. Since $x \in \widetilde {N(h)}$ we have $x \sim h$ and, $x$ and $h$ have same neighbors. Since $k$ is a neighbor of $h$ we have, $ k$ is a neighbor of $x$, i.e., $x \sim k.$  Combining these, we have either $\left<x\right>\subseteq \left<h\right>$ or $\left<h\right>\subseteq \left<x\right>$, and either $\left<x\right>\subseteq \left<k\right>$ or $\left<k\right>\subseteq \left<x\right>$. Suppose $\left<x\right>\subseteq \left<h\right>$. Then, either $\left<x\right>\subseteq \left<k\right>$ or $\left<k\right>\subseteq \left<x\right>$. If $\left<x\right>\subseteq \left<h\right>$ and $\left<k\right>\subseteq \left<x\right>$, this would imply $\left<k\right>\subset \left<h\right>$, contradicting our assumption that $\left<h\right>\subsetneq  \left<k\right>$. If $\left<x\right>\subseteq \left<h\right>$ and $\left<x\right>\subseteq \left<k\right>$ we have   $\left<x\right>\subseteq \left<h\right> \subseteq \left<k\right>,$   i.e., there exists a cyclic chain containing  $x, h$ and $ k.$   Next, suppose  $\left<h\right>\subseteq \left<x\right>$. Then, either $\left<x\right>\subseteq \left<k\right>$ or $\left<k\right>\subseteq \left<x\right>$. Together, this gives either $\left<h\right>\subseteq  \left<x\right> \subseteq \left<k\right>$ or $\left<h\right> \subseteq \left<k\right> \subseteq \left<x\right>$. Hence, in both the cases, $x,h$ and $k$ are contained in the some cyclic chain of subgroups of $G$.   
	\end{proof}    
\end{lemma}

\begin{lemma}\label{lab3}
	Let $G$ be a finite group and let  $\left<h_m\right>\subset \left<h_{m+1} \right> \subset \ldots \subset \left<h_{n-1} \right> \subset \left<h_n\right>$ be  a cyclic chain in $G$ where $\left<h_n\right>$ is a cyclic subgroup of prime-power order. Assume the following:
	\begin{enumerate}
		\item If an element $t\sim h_m$ or $h_n$ in $\mathscr{P}(G)$, then $t$  has a prime-power order.
		\item  If every maximal cyclic chain $\mathcal C$ of prime-power order subgroups containing $\left<h_m\right>$, then it also contains $\left<h_n\right>$. 
	  	\end{enumerate} Then, for any $i$ with $m \leq i \leq n$ we have,  $\bigcup\limits_{j=m}^{n}gen \left<h_j\right>\setminus \{h_j\}\subseteq \widetilde {N(h_i)}$.
  	\end{lemma}
	\begin{proof} Let $ t \in \bigcup\limits_{j=m}^{n}gen \left<h_j\right>\setminus \{h_j\}$ and let $i$ be such that $m \leq i \leq n.$   So, $t\in gen \left<h_k\right>$, i.e., $\left<t\right>= \left<h_k\right>$ for some $k$ with $m\leqslant k\leqslant n $. Note that  $\left<h_m\right>\subseteq \left<t\right>= \left<h_k\right>$ and $ \left<h_i\right> \subseteq \left<h_n\right> $. Since $\left<h_n\right>$ has prime-power order it has a unique maximal chain of subgroups descending from it. Hence, either  $\left<h_k\right> \subseteq \left<h_i\right>$ or $\left<h_i\right> \subseteq \left<h_k\right>.$
	We first assume that $\left<h_k\right> \subseteq \left<h_i\right>.$ The other case $\left<h_i\right> \subseteq \left<h_k\right>$ can be dealt with similarly.  
	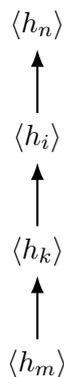
\begin{figure}[h!]
		\centering
		\begin{tikzpicture}[
		node distance=1.5cm,
		every node/.style={font=\small, align=center},
		arrow/.style={-Latex, thick}
		]
		\node (hm) {$\langle h_m \rangle$};
		\node (hk) [above of=hm] {$\langle h_k \rangle$};
		\node (hi) [above of=hk] {$\langle h_i \rangle$};
		\node (hn) [above of=hi] {$\langle h_n \rangle$};
		\draw[arrow] (hm) -- (hk);
		\draw[arrow] (hk) -- (hi);
		\draw[arrow] (hi) -- (hn);
		\end{tikzpicture}
		\caption{Chain of subgroups of prime-power order, showing inclusion relations.}
		\label{fig:chainC_simple}
	\end{figure}	
	We need to show, $t \in \widetilde {N(h_i)}$, i.e., $t \sim h_i$ and $N(h_i)\setminus \{t\}=N(t)\setminus \{h_i\}$ in $\mathscr{P}(G)$.
	Since  $\left<t\right>=\left<h_k\right> \subseteq \left<h_i\right>$, we have, $t \sim h_i$.	 Next, we prove that $N(h_i)\setminus \{t\}=N(t)\setminus \{h_i\}$ in $\mathscr{P}(G)$. To do this we first establish the following two claims.
	\noindent \underline{Claim 1}: Any maximal cyclic chain $\mathcal{C}$  of prime-power order subgroups that contains $\left<h_k\right>$ also contains $\left<h_i\right>$.\\
	 Let $\mathcal{C}$ be any maximal cyclic chain of subgroups of prime-power order with  $\left<h_k\right>\in \mathcal{C}$. We show $\left<h_i\right>\in \mathcal{C}$. Since  $ \left<h_m\right> \subseteq \left<h_k\right>  \in \mathcal C$ and  $\left<h_k\right>$ has prime-power order it follows that  $\left<h_m\right> \in \mathcal C.$ Since $\mathcal C$ is a maximal cyclic chain, hence by assumption (2) we obtain  $\left<h_n\right> \in \mathcal C$.  We thus have, $\mathcal{C}$ contains both $\left<h_m\right>$ and $\left<h_n\right>.$ Since  $\left<h_n\right>$ has prime-power order and $\mathcal C$ is maximal cyclic chain we obtain,   $\left<h_i\right> \in \mathcal C.$ This proves the claim 1. 
	 
	 \noindent \underline{Claim 2}: If $ x \sim h_i$ or $ x \sim t$ then, $x$ has a prime-power order. \\
	 If $ x \sim h_i$ then $\left<h_i\right> \subseteq \left<x\right> $ or $\left<x\right> \subseteq \left<h_i\right>.$ If $\left<h_i\right> \subseteq \left<x\right> $ then $\left<h_m\right> \subseteq \left<x\right>.$ This shows $x \sim h_m$ and hence by assumption (1),  $x$ has a prime-power order. If $\left<x\right> \subseteq \left<h_i\right>$ then $\left<x\right> \subseteq \left<h_n\right>.$ This gives $x \sim h_n$ and hence by assumption (1),  $x$ has a prime-power order. The case  $ x \sim t$ can be proved similarly. This proves Claim 2.

 We now show that $N(h_i)\setminus \{t\}=N(t)\setminus \{h_i\}$.  So, let $ x \in N(h_i) \setminus \{t\})$, i.e., $x\sim h_i$. Hence, $\left<x\right>\subseteq \left<h_i\right>$ or $\left<h_i\right>\subseteq \left<x\right>$. Suppose first that $\left<x\right>\subseteq \left<h_i\right>$. Already $\left<h_k\right>\subseteq \left<h_i\right>$ so we have, $x,h_k\in \left<h_i\right>$. Since $\left<h_i\right>$ is of order prime-power order we have either $\left<x \right> \subseteq \left<h_k\right>$ or  $\left<h_k \right> \subseteq \left< x\right>.$ Hence, $ x \sim h_k.$ Further, as $\left<h_k \right>  = \left<t \right> $ we obtain, $ x \sim t$, i.e., $x \in N(t) \setminus \{h_i\}.$ Next, we deal with the case $\left<h_i\right>\subseteq \left<x\right>.$ Since $\left<h_k\right>\subseteq \left<h_i\right>$ we have $\left<h_k \right>\subseteq \left<h_i \right> \subseteq \left<x\right>$ showing that $ x \sim h_k.$  Since $\left<h_k \right>  = \left<t \right> $ we obtain, $ x \sim t$, i.e., $x \in N(t) \setminus \{h_i\}.$ Hence, in both the cases we have, $x \in N(t) \setminus \{h_i\}$ proving that  $N(h_i)\setminus \{t\} \subseteq  N(t)\setminus \{h_i\}$. Thus, $N(h_i)\setminus \{t\}=N(t)\setminus \{h_i\}$. 
 
 We now prove the other containment. So, let $x \in N(t)\setminus \{h_i\}$, i.e., $ x \sim t$. Hence, $\left<x\right>\subseteq \left<t\right>= \left<h_k\right> $ or $\left<t\right>=\left<h_k\right>\subseteq \left<x\right>$.  If $\left<x\right>\subseteq \left<h_k\right> $, then since $\left<h_k\right>\subseteq \left<h_i\right> $ we have, $\left<x\right>\subseteq \left<h_i\right>.$ Thus, $x\in N(h_i) \setminus \{ t\}$. So, suppose now that $\left<t\right>= \left<h_k\right>\subseteq \left<x\right>$.  Note that by Claim 2 we find that $x$ is an element of prime-power order. So if  a maximal cyclic chain $\mathcal D$ of subgroups of prime-power order  contains  $\left<x\right>$ then $\mathcal D$ has to contain $\left<h_k\right>$. By claim 1, we obtain $\left<h_i\right> \in \mathcal D.$ Since $\mathcal D$ is a chain containing both $\left<h_i\right>$  and $\left<x\right>$ we find that either $\left<x\right> \subseteq \left<h_i\right>$ or $\left<h_i\right> \subseteq \left<x\right>.$ In any case we have, $ x \sim h_i$ showing that $ x \in N(h_i) \setminus \{ t\}$. Hence, $N(t)\setminus \{h_i\} \subseteq N(h_i)\setminus \{t\}.$ So, we have $N(h_i)\setminus \{t\}=N(t)\setminus \{h_i\}$ in $\mathscr{P}(G)$ and hence,  $ t \in \widetilde {N(h_i)}.$
             \end{proof}               

\begin{lemma}\label{lab4}
	Let $G$ be a finite group, and let $\left< e\right>  \neq \left<h\right>\subsetneq \left<k\right>\leqslant G$ be cyclic subgroups of prime-power order. Assume that any $t\sim h$ or $k$ in $\mathscr{P}(G)$, $t$ has prime-power order. Then:
	\begin{enumerate}
		\item If every maximal chain of cyclic subgroups of prime-power order containing $\left<h\right>$ also contains $\left<k\right>$ then, $N_h= N_k$.
		\item If there exists a maximal chain $\mathcal{C}$ of cyclic subgroups of prime-power order with $\left<h\right>\in \mathcal{C}$ and $\left<k\right>\notin \mathcal{C}$ then, $N_h\leq N_k$.
	\end{enumerate}

	\begin{proof} 
	 (1) We show that $\widetilde {N(h)}\setminus  \{k\} = \widetilde {N(k)}\setminus\{h\}.$ To do this, we first show that $\widetilde {N(h)}\setminus  \{k\} \subseteq  \widetilde {N(k)}\setminus\{h\}.$ So, let $x\in \widetilde {N(h) }\setminus  \{k\}.$ Since $x \sim h$ it follows from the hypothesis that $x$ has a prime power order. To show that $x \in \widetilde {N(k)}\setminus\{h\}$ we need to show $ x \sim k$ and $N(x) \setminus \{k\}= N(k) \setminus \{x\}.$ We first show that $ x \sim k.$ Since $ x \sim h$ we have $\left<x\right>  \subseteq  \left<h\right>$    or $\left<h\right>  \subseteq  \left<x\right>.$ If  $\left<x\right>  \subseteq  \left<h\right>$  then  $\left<x\right>  \subseteq  \left<h\right> \subseteq  \left<k\right>$ so it follows that $x \sim k.$ If  $\left<h\right>  \subseteq  \left<x\right>$ then we consider a maximal  chain of cyclic $\mathcal C$ subgroups of prime-power order  that contains both $\left<h\right> $ and $\left<x\right>.$  By hypothesis, chain $\mathcal C,$  as it contains $\left<h\right> $, it must contain $\left<k\right>.$ Hence, either $\left<x\right>  \subseteq  \left<k\right>$    or $\left<k\right>  \subseteq  \left<x\right>$ showing that $ x \sim k.$ We now show that $N(x) \setminus \{k\}= N(k) \setminus \{x\}.$ We first show that $ N(x) \setminus \{k\}\subseteq N(k) \setminus \{x\} .$ To do this, let $t \in N(x) \setminus \{k\}.$ Note that since $x\in \widetilde {N(h) }\setminus \{k\}$, by definition of $\widetilde {N(h) }$ we have that $x $ and $h$ have the same set of neighbors. So, $t \sim h$. By hypothesis $t$ has a prime-power order. So, consider a maximal chain $\mathcal D$ of cyclic subgroups  $\left<h\right> $ and $\left<t\right>$ both. From the hypothesis, $\left<k\right> \in \mathcal D.$ Thus,  $\left<t\right>  \subseteq  \left<k\right>$    or $\left<k\right>  \subseteq  \left<t\right>.$ This shows $ t \sim k$ proving that  $ N(x) \setminus \{k\}\subseteq N(k) \setminus \{x\}.$ Conversely, to show that $ N(k) \setminus \{x\} \subseteq N(x) \setminus \{k\}$ we let $ t \in N(k) \setminus \{x\}.$ By hypothesis, $t$ has prime-power order. We now claim that $ t \sim x.$ Since $ t \sim k$ we have $\left<k\right>  \subseteq  \left<t\right>$    or $\left<t\right>  \subseteq  \left<k\right>.$ If $\left<k\right>  \subseteq  \left<t\right>$ then $\left<h\right>  \subseteq  \left<t\right>$ giving $ h \sim t.$ As noted earlier $x$ and $h$ have the same set of neighbors, so $t \sim x.$ Suppose $\left<t\right>  \subseteq  \left<k\right>.$ In this case, since  $\left<k\right>$ has a prime-power order and both $\left<h\right>$, $\left<t\right>$ are subgroups of  $\left<k\right>$ we obtain that  $\left<h\right>  \subseteq  \left<t\right>$    or $\left<t\right>  \subseteq  \left<h\right>.$ Thus, $ t \sim h.$ As noted earlier $x$ and $h$ have the same set of neighbors, so $t \sim x.$ This shows that $ t \in N(x) \setminus \{k\}$ proving the containment $ N(k) \setminus \{x\} \subseteq N(x) \setminus \{k\}.$ So, $N(x) \setminus \{k\}= N(k) \setminus \{x\}.$ Hence, $x \in \widetilde {N(k)} \setminus \{h\}.$ This gives us the containment  $\widetilde {N(h)}\setminus  \{k\} \subseteq  \widetilde {N(k)}\setminus\{h\}.$ We now prove the other containment  $  \widetilde {N(k)}\setminus\{h\} \subseteq \widetilde {N(h)}\setminus  \{k\}.$ Let $ x \in  \widetilde {N(k)}\setminus\{h\}.$ To show $ x \in  \widetilde {N(h)}\setminus\{k\}$ it suffices to show that $ x \sim h$ and $ N(h) \setminus \{x\} = N(x) \setminus \{h\}.$  We first show that $ x \sim h.$ Note that since $ k \sim x$ we have $\left<k\right>  \subseteq  \left<x\right>$    or $\left<x\right>  \subseteq  \left<k\right>.$ If $\left<k\right>  \subseteq  \left<x\right>$ then $\left<h\right>  \subseteq  \left<x\right>$ and hence $ x \sim h.$ If $\left<x\right>  \subseteq  \left<k\right>$ then both $\left<x\right>$ and $\left<h\right>$ are subgroups of $\left<k\right>$ which is of prime-power order. Hence, $\left<x\right>  \subseteq  \left<h\right>$    or $\left<h\right>  \subseteq  \left<x\right>$ showing that $ x \sim h.$ We now show that $ N(h) \setminus \{x\} = N(x) \setminus \{h\}.$ To do this let $ t \in N(h) \setminus \{x\}.$ So,  $\left<t\right>  \subseteq  \left<h\right>$    or $\left<h\right>  \subseteq  \left<t\right>.$  By hypothesis, $t$ has a prime-power order.  So, let $\mathcal C$ be a maximal chain of order containing both $ \left<h\right>$  and  $\left<t\right>.$ By hypothesis, $\mathcal C$ must contain  $\left<k\right>.$  Hence,   $\left<k\right>  \subseteq  \left<t\right>$    or $\left<t\right>  \subseteq  \left<k\right>.$ So, $ t \sim k.$ Now, as $x\in \widetilde {N(k) }\setminus \{h\}$ we already have $N(k) \setminus \{x\} = N(x) \setminus \{k\}.$ So $t \sim x$ and $ N(h) \setminus \{x\} \subseteq N(x) \setminus \{h\}.$ To prove the reverse containment let $t \in N(x) \setminus \{h\}.$ So, $t \sim x$ and since $ N(x) \setminus  \{k \}= N(k) \setminus  \{x \}$  we obtain, $ t \sim k.$ So, $\left<k\right>  \subseteq  \left<t\right>$ or $\left<t\right>  \subseteq  \left<k\right>.$ If $\left<k\right>  \subseteq  \left<t\right>$ then $\left<h\right>  \subseteq  \left<t\right>$ giving us $t \sim h.$ If $\left<t\right>  \subseteq  \left<k\right>$   then we have both the subgroups $\left<t\right> $ and $\left<h\right> $contained in $\left<k\right> $. 
 Since $\left<k\right> $ is a subgroup of prime-power order we find that $\left<t\right>  \subseteq  \left<h\right>$ or $\left<h\right>  \subseteq  \left<t\right>.$ Hence, $t \sim h$, i.e., $t \in N(h) \setminus \{x\}. $ This proves the other containment $N(x) \setminus \{h\} \subseteq N(h) \setminus \{x\}.$ Thus,  $ N(h) \setminus \{x\} = N(x) \setminus \{h\}$ proving that  
 $\widetilde {N(k)}\setminus  \{h\} \subseteq  \widetilde {N(h)}\setminus\{k\}.$ So, $\widetilde {N(h)}\setminus  \{k\} = \widetilde {N(k)}\setminus\{h\}.$ Recall that $N_h= | \widetilde {N(h)}|+1. $  Notice that $ k \in   \widetilde {N(h)} $ if and only if $ h \in   \widetilde {N(k)}.$  This gives $N_h =N_k.$\\
 (2) Suppose there exist a maximal chain $\mathcal{C}$ of cyclic subgroups of prime-power order with $\left<h\right>\in \mathcal{C}$ and $  \left<k \right>\notin \mathcal{C}$. Suppose $|\left<h\right>|= p^{m_1}$ and $|\left<k\right>|= p^{m_2}$ with $m_1 < m_2$. From Lemma \ref {amar} we know that $  |gen (\left <k\right >)| \leq N_k,$ i.e.,  is $\phi(p^{m_2})=p^{m_2}- p^{m_2-1}= p^{m_2-1}(p-1) \leq N_k.$  We now claim that\\
 \underline{Claim:}   $N_h \leq p^{m_2-1}+1.$\\
 To prove this claim we first show that $\widetilde{ N(h)} \subseteq A=\left<k\right> \setminus gen(\left<k\right>).$  To prove the containment let $x \in \widetilde{ N(h)}. $ So, $ x \sim h$ and $ N(x) \setminus \{h\} =N(h) \setminus \{x\}.$ Since $ x \sim h$ we have  $\left<x\right>  \subseteq  \left<h\right>$  or $\left<h\right>  \subseteq  \left<x\right>.$ If $\left<x\right>  \subseteq  \left<h\right>$ then  $ x \in  \left<k\right>.$ Since $\left<h\right>\subsetneq \left<k\right>$  we see that $ x \notin \text{gen}(\left<k\right>)$ and hence $ x \in A$ in this case. Suppose $\left<h\right>  \subseteq  \left<x\right>.$ Since $ x \sim h$ it follows from the  hypothesis that $x$ is an element of prime-power order. We first show that $ x \notin gen(\left<k\right>).$ To see this suppose on the contrary that $ x \in gen(\left<k\right>).$ Note that as  $\left<k\right> \notin \mathcal C$ there exists $\left<y\right> \in \mathcal C$ such that $\left<y\right> \nsubseteq \left< x\right> $ and  $\left<x\right> \nsubseteq \left< y\right>.$ This gives $ y \nsim x.$  Since $\left<h\right> \subseteq \left< y\right> $ we have, $ y \sim h.$  So, $y \in N(h) \setminus \{x\} = N(x) \setminus \{h\} $ and hence we obtain, $ y \sim x. $ This contradicts  $ y \nsim x$ showing that $ x \notin gen(\left<k\right>)$. So, $ k \in N(h ) \setminus \{x\} = N(x ) \setminus \{h\}.$ Thus, $ k \sim x$ which gives us  $\left<x\right>  \subsetneq  \left<k\right>$ or $\left<k\right>  \subsetneq  \left<x\right>.$  We show that  $\left<k\right>  \subsetneq  \left<x\right>$ is not possible. To do this, we choose  as before $\left<y\right> \in \mathcal C$ such that $\left<y\right> \nsubseteq \left< x\right> $ and  $\left<x\right> \nsubseteq \left< y\right>.$ As before,  $y \in N(h) \setminus \{x\} = N(x) \setminus \{h\} $ and hence, $ y \sim x. $ This gives,   $\left<x\right>  \subseteq  \left<y\right>$ or $\left<y\right>  \subseteq  \left<x\right>$. If $\left<x\right>  \subseteq  \left<y\right>$ then we have, $\left<k\right>  \subseteq  \left<y\right>.$ Since $\mathcal C$ is a maximal chain of cyclic subgroups of prime-power order we find that $\left<k\right> \in \mathcal C$ which is a contradiction.  So,  $\left<x\right>  \subseteq  \left<y\right>$ is not possible. If $\left<y\right>  \subseteq  \left<x\right>$ then we have two distinct incomparable chains descending from $\left<x\right>$ namely,  $\left<h\right>  \subseteq  \left<y\right>  \subseteq  \left<x\right>$  and $\left<h\right>  \subseteq  \left<k\right>  \subseteq  \left<x\right>.$ Since $x$ has a prime-power order this is not possible. This shows that  $\left<k\right>  \subseteq  \left<x\right>$ is not possible.  This leaves us with  $\left<x\right>  \subsetneq  \left<k\right>$ giving us $x \in A.$ This shows that  $\widetilde{ N(h)} \subseteq A= \left<k\right> \setminus gen(\left<k\right>).$  Hence, $|\widetilde{ N(h)}| \leq |A|$ and so  $N_h=|\widetilde{ N(h)}|+1 \leq |A|+1.$ It is easy to see that $|A| = p^{m_2-1}.$ It now follows that  $N_h \leq p^{m_2-1}+1$ proving the claim. We now show that $N_h \leq N_k.$ We do this in two cases, namely, $p \geq 3$ and  $p=2.$ If $ p \geq 3,$ it is clear that $N_h \leq p^{m_2-1} +1 \leq p^{m_2-1}(p-1) \leq N_k.$ For $p=2,$ we argue as follows. We divide the case $ p=2$ into two subcases namely, $ e \in \widetilde{ N(h)}$ or $ e \notin \widetilde{ N(h)}.$ If  $ e \notin \widetilde{ N(h)}$ then  $\widetilde{ N(h)} \subseteq A \setminus \{e\}.$ So,  $|\widetilde{ N(h)}| \subseteq |A \setminus \{e\}|,$ giving us $N_h=|\widetilde{ N(h)}|+1 \subseteq |A \setminus \{e\}|+1=2^{m_2-1}-1+1=2^{m_2-1}.$ Already, $ 2^{m_2-1}(2-1)= 2^{m_2-1} \leq N_k.$ Thus, in this case we have, $N_h \leq N_k.$ Suppose now that $ e \in \widetilde{ N(h)}.$ Observe in this case that $ h \in U(\mathscr{P}(G)).$  Since  $ h \neq e$ we get that $|U(\mathscr{P}(G))| >1.$ If  $\mathscr{P}(G)$ is complete then by \cite[Theorem 2.12]{Chakrabarty 1} we have that $G$ is cyclic of prime of power order. In this case, it is immediate that $\widetilde{ N(h)} = N(h) \setminus \{h\}$, $\widetilde{ N(k)} = N(k) \setminus \{k\}$ and hence, $N_h=|G|=|N_k|.$ If $\mathscr{P}(G)$ is not a complete graph then by Proposition \ref{Cam} (2)(i), we have $\mathscr{P}(G)$ either $G$ is a cyclic group of non-prime-power order $n \geq 6$ or $G$ is a generalized quaternion group of order $2^n.$ When $G$ is a cyclic group of non-prime-power order $n \geq 6$,  $U(\mathscr{P}(G))$ consists of identity and generators of $G.$ Since $ h \in U(\mathscr{P}(G))$ and $h \neq e$ it follows that $h$ is a generator of $G.$ Note however that  $\left<h\right>\subsetneq \left<k\right>$ and hence $h$ cannot be generator of $G.$ Hence, $G$ is a generalized quaternion group of order $2^n.$ In this case, $U(\mathscr{P}(G))$ consists of $e$ and unique order two element in $G$ and so, $U(\mathscr{P}(G)) = \{e, h\}.$ Hence, $\widetilde{ N(h)}= \{e\}$ and consequently we have $N_h = |\widetilde{ N(h)}|+1=2.$ Note that since $\left<h\right>\subsetneq \left<k\right>$ we see that  order $k$ is at least $4.$ Hence,  $2=\phi(4) \leq |gen(\left<k\right>)|.$ By Lemma \ref{amar}, $2=|gen(\left<k\right>)| \leq N_k.$  So, $N_h =2 \leq |gen(\left<k\right>)| \leq N_k.$ This completes the proof. 
\end{proof}  
\end{lemma}
\begin{lem} \label{monotonic}
	Let $G$ be a finite non-cyclic group, and let $\left<h\right>\subseteq \left<k\right>\leqslant G$. Then, $N_h\leqslant N_k$.
	\begin{proof} 
		Case (i): $h$ and $k$ are elements of non-prime-power order.
		Suppose first that $h\neq e$. By Proposition \ref{non-prime-power}, we have $N_h= | gen (\left <h \right >)|$ and $N_k= | gen (\left <k \right >)|$. From Proposition \ref{cyclic}(1), we conclude that $N_h\leqslant N_k$.
		Now, suppose $h=e$. By Proposition \ref{Ne} if $G$ is a generalized quaternion group of order $2^n$ then, $N_h=2$. Otherwise, $N_h=1$. Since $k$ has nonprime power order it follows from Proposition \ref{non-prime-power} that  $N_k= | gen (\left <k \right >)|\geqslant 2$. Hence, $N_h\leqslant N_k$.\\
		Case (ii): The case where $k$ is element of  prime-power order and $h$ is element of non-prime-power order does not arise as $\left<h\right>\subseteq \left<k\right>.$ \\
		Case (iii): Suppose $h$ is an element of prime-power order and $k$ is element of non-prime-power. From Proposition \ref{Euler} it follows that, $N_h= | gen (\left <h \right >)|$ and from  Proposition  \ref{non-prime-power} we obtain, $N_k= | gen (\left <k \right >)|$. Now since  $\left<h\right>\subseteq \left<k\right>$  it follows from Proposition \ref{cyclic} that $N_h\leqslant N_k$. \\
		Case (iv): Suppose $h$ and $k$ are elements of order prime-power, i.e., $ |\left<h\right>|= p^m$ and  $ |\left<k\right>|= p^n$ where $m\leqslant n$ and $p$ is a prime.\\
		Subcase (i): Suppose there exists $l \in G$ such that $\left<k\right>\subset \left<l\right>$ and $l$ is an element of non-prime-power order.
		Since  $\left<h\right> \subset \left<k\right>\subset \left<l\right>$ by Proposition \ref{Euler} , $N_h= | gen (\left <h \right >)|$ and $N_k= | gen (\left <k \right >)|$. Hence, from Proposition  \ref{cyclic} it follows that, $N_h\leqslant N_k$.\\
		Subcase (ii): Suppose there does not exist any $l \in G$ such that $\left<k\right>\subset \left<l\right>$ and $l$ is an element of non-prime-power order. By Proposition \ref{amar}, $| gen (\left <k \right >)| \leq N_k.$ Suppose there exists $h' \in G$ such that $\left<h\right>\subset \left<h'\right>$ and $h'$ is an element of non-prime-power order. In this case,  by Proposition \ref{Euler} we get $N_h= | gen (\left <h \right >)|.$ So, by Proposition \ref{cyclic} we find that $N_h= | gen (\left <h \right >)| \leq | gen (\left <k \right >)| = N_k.$  Now suppose that there does not exist $h'$ with $\left<h\right>\subset \left<h'\right>$ and $h'$ an element of non-prime-power order. Now here we have $\left<h\right>\subseteq \left<k\right>\leqslant G$, where both $\left<h\right>$ and $\left<k\right>$ are cyclic subgroups of prime-power order.
		
		 Since there is no cyclic subgroup of non-prime-power order containing $\left<h\right>$ (and hence containing $k$) we have that if an  element $t\sim h$ or $k$ in the $\mathscr{P}(G)$, then $t$ has a prime-power order. We now have two cases as follows.\\
		  Case (A): Any maximal chain $\mathcal C$ of cyclic subgroups of prime-power order containing $\left<h\right>$ also contains $\left<k\right>$. \\		   
		   Case (B): There exists a maximal chain $\mathcal C$ of cyclic subgroups of prime-power order containing $\left<h\right>$  but not containing $\left<k\right>.$\\		   
		   In the case (A) we have, by Lemma \ref{lab4}(i), $N_h=N_k$. In the case (B) it follows from Lemma \ref{lab4}(ii) that $N_h \leqslant N_k$. Thus, in both the cases we conclude that $N_h\leqslant N_k$. Hence the proof is complete.                                                    
	\end{proof} 
\end{lem}

\section{Main Results}

\begin{lemma}\label{ordera,b}
Let $G$ be any finite group. If $a, b \in \left<c\right>$ and $|\langle a\rangle|$ divides $|\langle b\rangle|$ or $|\langle b\rangle|$ divides $|\langle a\rangle|$, then $a \sim b $ in $\mathscr{P}(G)$.
\end{lemma}
\begin{proof}
Suppose  $a, b \in \left<c\right>,$ i.e., $\langle a\rangle, \langle b\rangle \subseteq \left<c\right>$. Set  $d_1=|\langle a\rangle|$ and  $d_2=|\langle b\rangle|.$ Since $\left<c\right>$ is cyclic, $\langle a \rangle$ and $\langle b \rangle$ are unique  cyclic subgroups of $\left<c\right>$ of  orders $d_1$ and $d_2$ respectively. Without loss of generality assume that $d_1$ divides $d_2.$ Since $\langle b\rangle$ is cyclic of order $d_2$ and $d_1 | d_2$ we obtain a unique subgroup $\langle h\rangle$  of order $d_1$ inside $\langle b\rangle.$ Thus, $\langle a\rangle$ and $\langle h\rangle$ are subgroups of order $d_1$ inside $\langle c \rangle.$  Hence, by uniqueness of subgroups, $\langle a\rangle = \langle h\rangle.$ This shows  $\langle a\rangle \subseteq \left<b\right>$ and so, $a \sim b $ in $\mathscr{P}(G)$. 
\end{proof}

\begin{theorem}\label{main thm 1}
Let $G$ be a finite non-cyclic group. For $a, b \in G, $ assume that  $ a \nsim b $ in $\mathscr{P}(G)$ and $N_a \neq N_b.$ The following statements are equivalent.
 \begin{enumerate}
\item $a \sim b $ in $\mathscr{P}_{e}(G)$  
\item $ \exists ~ c  \in G$ such that $a, b \sim c$ in $\mathscr{P}(G)$ and $N_c \geqslant \max\{ N_a, N_b\}$. 
\end{enumerate}
\begin{proof}
Suppose $a \sim b $ in $\mathscr{P}_{e}(G)$. By definition of $\mathscr{P}_{e}(G)$, $ \exists ~ c  \in G$ such that $a,b \in \left<  c \right >$. Therefore, $ \left<  a \right >, \left<  b \right > \subseteq  \left<  c \right > $. By Lemma \ref{monotonic}, $ N_a, N_b \leqslant N_c $. So, it follows that  $N_c \geqslant \max\{ N_a, N_b\}$. Conversly, suppose that  $ \exists ~ c \in G$ such that $a, b \sim c$ in $\mathscr{P}(G)$ and $N_c \geqslant \max\{ N_a, N_b\}$. We claim that $a \sim b $ in $\mathscr{P}_{e}(G)$. 
By definition of $\mathscr{P}(G)$, as $a \sim c $ we have, $\left<a\right>\subseteq \left<c\right>$ or $\left<c\right>\subseteq \left<a\right>$. Similarly, since $b \sim c $ we have, $\left<b\right>\subseteq \left<c\right>$ or $\left<c\right>\subseteq \left<b\right>$. There are thus four cases as follows: \\
Case(i): $\left<a\right>\subseteq \left<c\right>$ and $\left<b\right>\subseteq \left<c\right>$, \hspace{1.2 cm}
Case(ii): $\left<a\right>\subseteq \left<c\right>$ and $\left<c\right>\subseteq \left<b\right>$\\
Case(iii): $\left<c\right>\subseteq \left<a\right>$ and $\left<b\right>\subseteq \left<c\right>$, \hspace{1 cm}
Case(iv): $\left<c\right>\subseteq \left<a\right>$ and $\left<c\right>\subseteq \left<b\right>$\\
Note that in  the cases (ii) and (iii)  we find that $a \sim b$ in  $\mathscr{P}(G)$ contradicting the assumption that $ a \nsim b$ in $\mathscr{P}(G)$. Hence, we are in either Case(i) or Case(iv), i.e., either $\left<a\right>, \left<b\right>\subseteq \left<c\right>$ or $\left<c\right> \subseteq \left<a\right>, \left<b\right>$. If $\left<c\right> \subseteq \left<a\right>, \left<b\right>$ then by Lemma \ref{monotonic}, $N_c\leq  N_a$ and $N_c\leq  N_b$. Since $N_a \neq N_b$ we obtain $N_c < \max\{ N_a, N_b\}$  contradicting the assumption that $N_c \geq \max\{ N_a, N_b\}$. Hence, $\left<a\right>, \left<b\right>\subseteq \left<c\right>$. It is now clear that  $a \sim b $ in $\mathscr{P}_{e}(G)$. This completes the proof.                    
\end{proof}
\end{theorem}

\begin{theorem}\label{main thm 2}
Let $G$ be a finite non-cyclic group. For $a, b \in G$ assume that $ a \nsim b $ in $\mathscr{P}(G)$ and $N_a = N_b= k.$ The following statements are equivalent. \begin{enumerate}
	\item $a \sim b $ in $\mathscr{P}_{e}(G)$  
	\item $\exists ~ c $ in $G$ such that $a, b\sim c$ in $\mathscr{P}(G)$ and $ N_c > k$.
\end{enumerate}
\begin{proof}
Suppose that, $a \sim b $ in $\mathscr{P}_{e}(G)$. So,  by definition of $\mathscr{P}_{e}(G)$, there exists $ c \in G $  such that $\left<a\right>, \left<b\right>\subseteq \left<c\right>$. From Lemma \ref{monotonic}, it follows that $N_a\leqslant N_c$ and  $N_b\leqslant N_c$. Since $N_a = N_b= k $, this gives $k\leqslant N_c,$ i.e., either $k< N_c$ or $k=N_c$.  We show that $k=N_c$ is not possible. So, assume that  $k=N_c$, i.e., $N_a=N_b=k=N_c$. Suppose first that $c$ is an element of prime-power order. Since $\left<a\right>, \left<b\right>\subseteq \left<c\right>$, we obtain that either  $\left<a\right> \subseteq  \left<b\right>$ or $\left<b\right> \subseteq  \left<a\right>.$ Hence, $ a \sim b$ in $\mathscr{P}(G)$ contradicting our assumption that $a \nsim b $ in $\mathscr{P}(G)$. Therefore, $c $ is an element of non prime-power order. If $a$ is an element of prime-power order then it follows from Proposition \ref{Euler} that $N_a= | gen (\left <a \right >)|,$ while if  $a$ is  an element of non prime-power order then we see from Proposition \ref{non-prime-power} that $N_a= | gen (\left <a \right >)|.$ Similarly,  $N_b= | gen (\left <b \right >)|.$ Further, it follows from Proposition \ref{non-prime-power} that  $N_c= | gen (\left <c \right >)|$. We thus have $\left<a\right>, \left<b\right>\subseteq \left<c\right>$ and $ | gen (\left <a \right >)|= | gen (\left <b \right >)|= | gen (\left <c \right >)|$.
Since  $\left<a\right> \subseteq \left<c\right>$ and  $ | gen (\left <a \right >)|= | gen (\left <c \right >)|,$ Proposition \ref{cyclic}  gives either $|\left<a\right>| = |\left<c\right>|$ or $|\left<c\right>| = 2|\left<a\right>|.$ Similarly, $\left<b\right> \subseteq \left<c\right>$ and  $ | gen (\left <b \right >)|= | gen (\left <c \right >)|$ gives either $|\left<b\right>| = |\left<c\right>|$ or $|\left<c\right>| = 2|\left<b\right>|.$ There are thus four cases and we settle them as follows.\\
\underline{Case (i)}: $|\left<a\right>| = |\left<b\right>|= |\left<c\right>|.$  In this case we have, $\left<a\right> = \left<b\right>= \left<c\right>$ giving us $ a \sim b$ in $\mathscr{P}(G).$ This contradicts our assumption that $ a \nsim b$ in $\mathscr{P}(G).$ \\
\underline{Case (ii)}: $|\left<a\right>| = 2 |\left<b\right>|= |\left<c\right>|.$ In this case, $|\left<b\right>|$ divides $|\left<a\right>|$ and so by Lemma \ref{ordera,b} we get $ a \sim b$ in  $\mathscr{P}(G),$ a contradiction. \\
\underline{Case (iii)}: $|\left<c\right>| = 2 |\left<a\right>|= |\left<b\right>|.$ In this case, $|\left<a\right>|$ divides $|\left<b\right>|$ and so by Lemma \ref{ordera,b} we get $ a \sim b$ in  $\mathscr{P}(G),$ a contradiction.\\
\underline{Case (iv)}: $|\left<c\right>| = 2 |\left<a\right>|= 2|\left<b\right>|.$ Clearly, $|\left<a\right>|= |\left<b\right>|$ and as $\left<a\right>, \left<b\right> \subseteq \left<c\right>$ we have, $\left<a\right> =\left<b\right>.$  So, we obtain  $ a \sim b$ in  $\mathscr{P}(G),$ a contradiction.\\
 We thus see that the case $ k=N_c$ is not possible and hence we have, $ k < N_c$ as desired.
 \newline \indent Next, suppose there exists $c \in G$ such that $a, b \sim c$ in $\mathscr{P}(G)$ and $ N_c > k$. Since $a \sim c $ in $\mathscr{P}(G)$, either $\left <a\right > \subseteq \left <c\right > $ or $\left <c\right > \subseteq \left <a\right >.$  Similarly, as $b\sim c $ we have, either $\left <b\right > \subseteq \left <c\right > $ or $\left <c\right > \subseteq \left <b\right >.$  We thus have four cases as in Theorem \ref{main thm 1} two of which can be ruled out as in the proof of Theorem \ref{main thm 1}. Hence, we are left with either $\left<a\right>, \left<b\right>\subseteq \left<c\right>$ or $\left<c\right> \subseteq \left<a\right>, \left<b\right>$. Suppose $\left<c\right> \subseteq \left <a\right >, \left <b\right >$. In this case we obtain from Lemma \ref{monotonic} that,  $ N_c\leq N_a$ and $ N_b \leq N_c$  giving us $ N_c\leq k $. This contradicts our assumption that $ N_c > k$. Hence, $\left<a\right>,\left<b\right> \subseteq \left<c\right>$ giving us $a\sim b$ in $\mathscr{P}_{e}(G)$ as desired. This completes the proof.                 
\end{proof}
\end{theorem}

\section{Algorithm for constructing $\mathscr{P}_{e}(G)$  from $\mathscr{P}(G)$  }
In this section, we give an algorithm for constructing the enhanced power graph  $\mathscr{P}_{e}(G)$  from the power graph $\mathscr{P}(G)$ without the knowledge of the underlying finite group $G.$ We then illustrate our algorithm with the help of an example. Recall that both the graphs $\mathscr{P}(G)$ and $\mathscr{P}_{e}(G)$ have the same vertex set $G.$ Further, $\mathscr{P}(G)$ is a subgraph of $\mathscr{P}_{e}(G).$ Hence, if $a \nsim b$ in $\mathscr{P}(G)$ then we need to decide  whether or not  $ a \sim b $ in $\mathscr{P}_{e}(G)$. \\

\begin{algorithm}
	Construction of the enhanced power graph $\mathscr{P}_{e}(G)$ directly  from the power graph $\mathscr{P}(G)$ for any finite group $G.$  Suppose that a graph $X$ is given such that $X \cong \mathscr{P}(G)$ for some finite group $G.$ Our aim to construct graph  $Y$ such that $Y \cong \mathscr{P}_e(G).$ In what follows, we denote by   $U(X)$ the set of vertices of $X$ which are joined to all other vertices. 
	\end{algorithm}

	\noindent {\bf Step (1)}: If $X$ is complete  then $X$ is the graph of a cyclic group of prime-power order by \cite[Theorem 2.12]{Chakrabarty 1}. Consequently, $Y \cong \mathscr{P}_{e}(G)$ is also complete and the algorithm terminates.\\
	
	\noindent {\bf Step (2)}: If $X$ is not a complete graph and $|U(X)|>2,$ then  $X$ is the graph $\mathscr{P}(G)$ of a cyclic group $G$ of non-prime-power order by Proposition \ref{Cam}(2)(i). Hence, it follows from \cite[Theorem 2.4]{Bera} that $Y \cong \mathscr{P}_e(G)$ is a complete graph. \\
	
	\noindent {\bf Step (3)}: If $X$ is not a complete graph and $|U(X)| \leq 2$, then, $X$ is the graph $\mathscr{P}(G)$ of non-cyclic group $G$ by Proposition \ref{Cam}. \\
	
	\noindent {\bf Step (4)}: For $ a, b \in X$ if $ a \sim b$ in $X$ then $ a \sim b$ in $Y \cong \mathscr{P}_e(G).$ If $ a \nsim b$ in $X$  then we proceed to the following steps to decide if $a$ is adjacent to $b$ in the enhanced power graph $ Y \cong \mathscr{P}_e(G).$\\
	
	\noindent {\bf Step (5)}: For each element $a \in G$, calculate the numbers $N_a$.  \\

	\noindent {\bf Step (6)}: Suppose $ a \nsim b$ in $X$ and $N_a= N_b=k.$ In this case, if there exists $ c \in X$ such that $a, b \sim c$ in $X$ with $ N_c > k$ then by Theorem \ref{main thm 2} we have, $ a \sim b$ in $ Y \cong \mathscr{P}_e(G).$ If such an element $c \in X$ does not exist then $ a \nsim b$ in the enhanced power graph.\\
	
	\noindent {\bf Step (7)}:  Suppose $ a \nsim b$ in $X$  and $N_a \neq N_b.$ In this case, if there exists $ c \in X$ such that $a, b \sim c$ in $X$ with $N_c \geqslant \max\{ N_a, N_b\}$ then by Theorem \ref{main thm 1} we have, $ a \sim b$ in $ Y \cong \mathscr{P}_e(G).$ If such an element $c \in X$ does not exist then $ a \nsim b$ in the enhanced power graph.\\

\noindent  \textbf{Illustration of the Algorithm:} We illustrate our algorithm for the  power graph $X$ in the figure below associated to  group $G$ having $12$ vertices.\\

\begin{center}
	
	\begin{figure}[h!]  \label{figure}
	\centering
	\begin{tikzpicture}[scale=1.6, line cap=round, line join=round]
		\tikzstyle{vertex}=[circle,fill=black,minimum size=7pt,inner sep=0pt]
		\foreach \i in {1,...,12}{
			\node[vertex] (v\i) at ({360/12*(\i-1)}:1.8) {};
			\node at ({360/12*(\i-1)}:2.1) {$v_{\i}$};
		}
		\foreach \j in {2,...,12}{
			\draw (v1)--(v\j);
		}
		\draw (v9)--(v10);
		\draw (v9)--(v11);
		\draw (v10)--(v12);
		\draw (v11)--(v12);
		\draw (v10)--(v11);
		\draw (v9)--(v12);
		\draw (v8)--(v11);
		\draw (v8)--(v12);
	     \draw[dashed] (v8)--(v9);
		\draw[dashed] (v8)--(v10);
		
	\end{tikzpicture}
	\caption{Power graph  and the enhanced power graph for dihedral group $D_6$ with additinal edges of $\mathscr{P}_e(G)$  indicated by dotted lines.}
\end{figure}
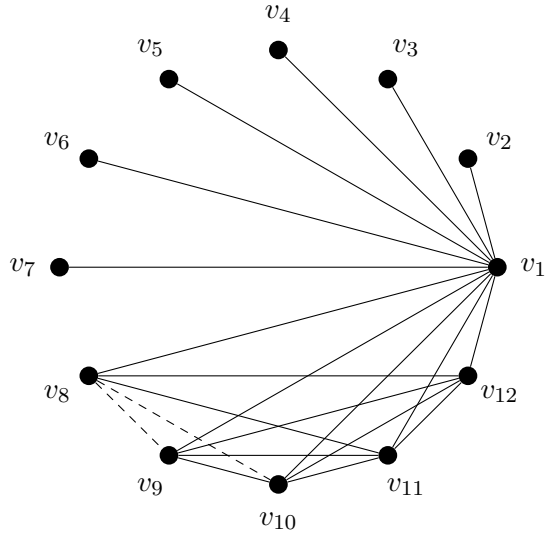

\end{center}
 To do this, we first determine whether he underlying  group $G$ is cyclic or non-cyclic for the given graph $X$. As $X$ is not complete, we conclude that underlying group $G$ is not a cyclic group of prime-power order. So, we proceed to find $|U(X)|$. In $X$, vertex $v_1$ is the only vertex that is adjacent to all other elements of $G$ and so, $|U(X))|= 1.$ Since $X$ is not complete and $|U(X))|= 1 < 2$, by step (3), the underlying group  $G$ is non-cyclic.  From step (4), for all $i\neq j$, $ v_i, v_j \in X:$\\ 
(i) if $ v_i \sim v_j$ in $X$, then $ v_i \sim v_j$ in $Y \cong \mathscr{P}_e(G);$ \\
(ii) if $ v_i \nsim v_j$ in $X$, then we proceed to step (5) to decide whether $v_i$ is adjacent to $v_j$ in the enhanced power graph $ Y \cong \mathscr{P}_e(G).$\\
By step (5), for each vertex $v_i$, we determine the number $N_{v_i}$. We have: $ N_{v_1}= 1, N_{v_2}= N_{v_3} = N_{v_4} = N_{v_5} = N_{v_6} = N_{v_7} = N_{v_8} =1$, $N_{v_9}=N_{v_{10}}= 2$ and $N_{v_{11}}=N_{v_{12}}= 2$. Vertices $v_8$ and $v_9$ are not adjacent in $X$. However, $N_{v_8}= 1,  N_{v_9}= 2$ and there exists a vertex $v_{11}$ such that $v_8$ and $v_9$ are both adjacent to $v_{11}$ in $X$, with $N_{v_{11}}= 2$. Thus, from step (6) we have, $v_8 \sim v_9$ in $Y$. Similarly, $v_8 \sim v_{10}$ in $Y$. It is easy to see that no additional edges are added to $ Y\cong \mathscr{P}_e(G)$ beyond these.\\

\noindent \textbf{Algorithm for Constructing the Difference Graph:}

\noindent Recall that the difference graph $\mathscr{D}(G)$ is defined to be the graph $\mathscr{P}_{e}(G) \setminus \mathscr{P}(G)$ with isolated vertices removed. The notion of difference graph of a group was proposed in \cite[Section 3.2]{Cameron 1} and formally introduced in \cite[Definition 1.3]{Cameron 5}. It is clear from the results of the previous section that diiference graph can be constructed from the power graph without the knowledge of the underlying group $G.$ We formally state this in the form of an algorithm as follows. \\

\begin{algorithm}
	 Construction of the Difference graph $\mathscr{D}(G)$  from  the power graph $\mathscr{P}(G)$  for any finite group $G$.
	 
	 \noindent \textbf{Step (1):} Use the algorithm given in Section 3 to construct $\mathscr{P}_{e}(G)$ from $\mathscr{P}(G)$.\\
	 \noindent \textbf{Step (2):} Consider the 
	 the graph $\mathscr{P}_{e}(G) \setminus \mathscr{P}(G).$ \\
	 \noindent \textbf{Step (3):} Remove the isolated vertices, if any, from the graph $\mathscr{P}_{e}(G) \setminus \mathscr{P}(G).$   The resulting graph is the difference graph $\mathscr{D}(G).$      

\end{algorithm}


\begin{thebibliography}{6}
	
\bibitem{Cameron 4} {G. Aalipour, S. Akbari, P. J. Cameron, R. Nikandish, and F. Shaveisi, On the structure of the power graph and the enhanced power graph of a group, \textit{Electron. J. Comb.} {\bf 24} (2017), P3.16.}

\bibitem{Bera} {S. Bera and A. K. Bhuniya, On enhanced power graphs of finite groups, \textit{J. Algebra Appl.} {\bf 17} (2018), 1850146.}

\bibitem{Bubboloni}{D. Bubboloni and N. Pinzauti, Critical classes of power graphs and reconstruction of directed power graphs,  \textit{J. Group Theory} {\bf 28} (2024), 713-739.}

\bibitem{Cameron 5}{S. Biswas, P. J. Cameron, A. Das and H. K. Dey, On the difference of the enhanced power graph and the power graph of a finite group, \textit{J. Combin. Theory Series A} {\bf 208} (2024), 105932.}
	
\bibitem{Brauer} {R. Brauer and K.A. Fowler, On groups of even order, \textit{The Annals of Mathematics} {\bf 62}(3) (1955), 567-583.}
	
\bibitem{Cameron 2} {P. J. Cameron, The power graph of a finite group II, \textit{J. Group Theory} {\bf 13} (2010), 779-783.}

\bibitem{Cameron 3} {P. J. Cameron and S. Ghosh, The power graph of a finite	group, \textit{Discrete Math.} {\bf 311} (2011), 1220-1222.}


\bibitem{Cameron 6} {P. J. Cameron and N. Maslova, Criterion of unrecognizability of a fnite group by its Gruenberg-Kegel graph, \textit{J. Algebra} {\bf 607} (2022), 186-213.}

\bibitem{Cameron 1} {P. J. Cameron, Graphs defined on groups, \textit{Int. J. Group Theory} {\bf 11} (2022), 53-107.}

\bibitem{Cameron 6} {P. J. Cameron, 	What can graphs and algebraic structures say to each other?, \textit{AKCE Int. J. Graphs Combin.} {\bf 21}(3) (2024), 249-254.}

\bibitem{Chakrabarty 1} {I. Chakrabarty, S. Ghosh and M. K. Sen, Undirected power graphs of semigroups, \textit{Semigroup Forum} {\bf 78} (2009), 410-426.}




\bibitem{Das}{B. Das, J. Ghosh and A. Kumar, The isomorphism problem of power graphs and a question of Cameron, \textit{$44^{th}$ IARCS  Annual Conference on Foundations of Software Technology and Theoretical Computer Science
		(FSTTCS)}(2024), https://arxiv.org/pdf/2305.18936}


\bibitem{Kelarev1}{A. V. Kelarev and S. J. Quinn, A combinatorial property and power graphs of groups, \textit{Contrib. General Algebra} {\bf 12} (2000), 229-235.}

\bibitem{Kelarev2}{A. V. Kelarev and  S. J. Quinn, Directed graph and combinatorial properties of semigroups, \textit{J. Algebra} {\bf 251} (2002), 16-26.}





\bibitem{Cameron6}{A. Kumar, L. Selvaganesh, P. J. Cameron and T. Tamizh Chelvam, Recent developments on the power graph of finite groups – a survey,  \textit{ AKCE Int. J. Graphs Combin.} {\bf 18} (2021), 65-94.}



\bibitem{Zahi}{S. Zahirović, I. Bošnjak, and R. Madarász, A study of enhanced power graphs of finite groups,  \textit{ J. Algebra Appl.} {\bf 19} (2020), 2050062.}


\end{thebibliography}
\end{document}